%% file: Zconvex.tex
\documentclass{amsart} 

\usepackage{a4wide}
\usepackage{graphicx}
\usepackage{color}

\newtheorem{theorem}{Theorem}[section]

\newtheorem{property}{Property}[section]
\newtheorem{lemma}[theorem]{Lemma}
\theoremstyle{definition}

\theoremstyle{remark}

\numberwithin{equation}{section}

\author[E. Duchi]{Enrica Duchi}
\address{Universit\'e Paris 7}
\email{duchi@liafa.jussieu.fr}
\thanks{The authors acknowledge support from the french ANR under
the SADA project.}

\author[S. Rinaldi]{Simone Rinaldi}
\address{Universit\`a di Siena}
\email{rinaldi@unisi.it}

\author[G. Schaeffer]{Gilles Schaeffer}
\address{CNRS, \'Ecole Polytechnique}
\email{schaeffe@lix.polytechnique.fr}

\title[\textsf{Z}-convex Polyominoes]{The number of \textsf{Z}-convex
polyominoes}



\begin{document}

\begin{abstract}
In this paper we consider a restricted class of polyominoes that we
call \textsf{Z}-convex polyominoes. \textsf{Z}-convex polyominoes are
polyominoes such that any two pairs of cells can be connected by a
monotone path making at most two turns (like the letter
\textsf{Z}). In particular they are convex polyominoes, but they
appear to resist standard decompositions. We propose a construction by
``inflation'' that allows to write a system of functional equations
for their generating functions. The generating function $P(t)$ of
\textsf{Z}-convex polyominoes with respect to the semi-perimeter turns
out to be algebraic all the same and surprisingly, like the generating
function of convex polyominoes, it can be expressed as a rational
function of $t$ and the generating function of Catalan numbers.


\end{abstract}

\maketitle 

\section{Introduction}

\subsection{Convex polyominoes}
In the plane $\Bbb Z \times \Bbb Z$ a {\em cell} is a unit square, and
a {\em polyomino} is a finite connected union of cells having no cut
point. Polyominoes are defined up to translations. A {\em column}
({\em row}) of a polyomino is the intersection between the polyomino
and an infinite strip of cells lying on a vertical
(horizontal) line. For the
main definitions and results concerning polyominoes we refer to
\cite{stan} and, for french aware readers, to \cite{mbm}.
Invented by Golomb \cite{dbintr39} who coined the term {\em
polyomino}, these combinatorial objects are related to many
mathematical problems, such as tilings~\cite{Gi, Go}, or games~\cite{Ga}
among many others. The enumeration problem for general polyominoes is
difficult to solve and still open. The number $a_n$ of polyominoes
with $n$ cells is known up to $n=56$ \cite{jensen-guttmann} and asymptotically,
these numbers satisfy the relation $\smash{\lim_n \left(a_n
\right)^{1/n}=\mu}$, \ $3.96 < \mu <4.64$, where the lower bound is a
recent improvement of \cite{bmrr}.

\begin{figure}[h]
\centerline{\hbox{\includegraphics[width=3in]{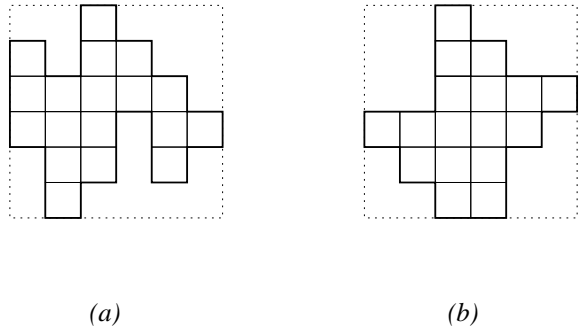}}}
 \caption{$(a)$ a column-convex (but not convex)
polyomino; $(b)$ a convex polyomino.}\label{convi}
\end{figure}

In order to probe further, several subclasses of polyominoes have been
introduced on which to hone enumeration techniques. One natural
subclass is that of {\em convex} polyominoes. A polyomino is said to
be {\em column-convex} [{\em row-convex}] when its intersection with
any vertical [horizontal] line of cells in the square lattice is
connected (see Fig.~\ref{convi} (a)), and {\em convex} when it is both
column and row-convex (see Fig.~\ref{convi} (b)). The {\em area} of a
polyomino is just the number of cells it contains, while its {\em
semi-perimeter} is half the length of the boundary. Thus, in a convex
polyomino the semi-perimeter is the sum of the numbers of its rows and
columns. Moreover, any convex polyomino is contained in a rectangle in
the square lattice which has the same semi-perimeter (called the {\em
minimal bounding rectangle} of the polyomino).

\medskip

The number $f_{n}$ of convex polyominoes with semi-perimeter $n+2$ was
obtained by Delest and Viennot, in \cite{DV}:
$$
f_{n+2}=(2n+11)4^n -4(2n+1) {2n \choose n}, \quad n \geq 0; \quad
f_{0}=1, \quad f_{1}=2.
$$
In particular the generating function of convex polyominoes with
respect to the semi-perimeter
$$F(t)=\sum_{n\geq0}f_{n}t^{n+2} = t^2+2t^3+7t^4+28t^5+120t^6+528t^7+O(t^8)$$
is an algebraic series which is has a rational expression
$R_0(t,d(t))$ in $t$ and the Catalan generating function
\[
c(t)=\frac{1-\sqrt{1-4t}}{2t}=1+t+2t^2+5t^3+14t^4+42t^5+132t^6+O(t^7).
\]
More precisely, the generating function of convex polyominoes with
respect to the numbers of columns (variable $x$) and rows (variable
$y$) is 
\[
F(x,y)=
\frac{8x^2y^2d(x,y)}{\Delta^2}+\frac{xy(1-x-xy-y)}{\Delta},
\]
where $d(x,y)$ is the unique power series satisfying the relation
$d=(x+d)(y+d)$,
\[
d(x,y)=\frac12(1-x-y-\sqrt{\Delta}),
\]
and,
\[
\Delta=(1-x-y)^2-4xy=(1-x-y)^2(1-\textstyle\frac{4xy}{(1-x-y)^2}).
\]
Observe that $d(t):=d(t,t)$ is just a shifted version of the Catalan
generating function, 
\[
d(t)=t(c(t)-1)=\frac12(1-2t-\sqrt{1-4t})
=t^2+2t^3+5t^4+14t^5+42t^6+132t^7+O(t^8).
\]
Incidentally, $d(x,y)$ is the generating function of parallel
polyominoes with respect to the numbers of columns and rows. 


\subsection{Monotone paths and $k$-convexity}
\begin{figure}
\centerline{\hbox{\includegraphics{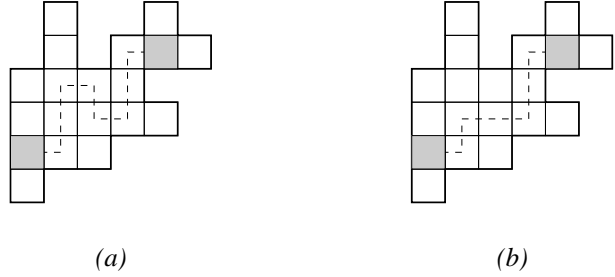}}} \caption{\label{paths}
$(a)$ a path between two highlighted cells in a polyomino; $(b)$ a
monotone path between the same cells, made only of north and east
steps.}
\end{figure}

In~\cite{lconv} the authors observed that convex polyominoes have the
property that every pair of cells is connected by a \emph{monotone
path}.  More precisely, a {\em path} in a polyomino is a self-avoiding
sequence of unitary steps of four types: north $N=(0,1)$, south
$S=(0,-1)$, east $E=(1,0)$, and west $W=(-1,0)$. A path is {\em
monotone} if it is made with steps of only two types.  Given a path
$w=u_1 \ldots u_k$, with $u_i \in \{N, S, E, W \}$, each pair of steps
$u_iu_{i+1}$ such that $u_i \neq u_{i+1}$, $0< i < k$, is called a
{\em change of direction}. These definitions are illustrated by
Fig.~\ref{paths}, in which the non monotone path (a) has 6 changes of
direction and the monotone path (b) has 4 changes of direction.  

\begin{figure}
\centerline{\hbox{\includegraphics[width=4.5in]{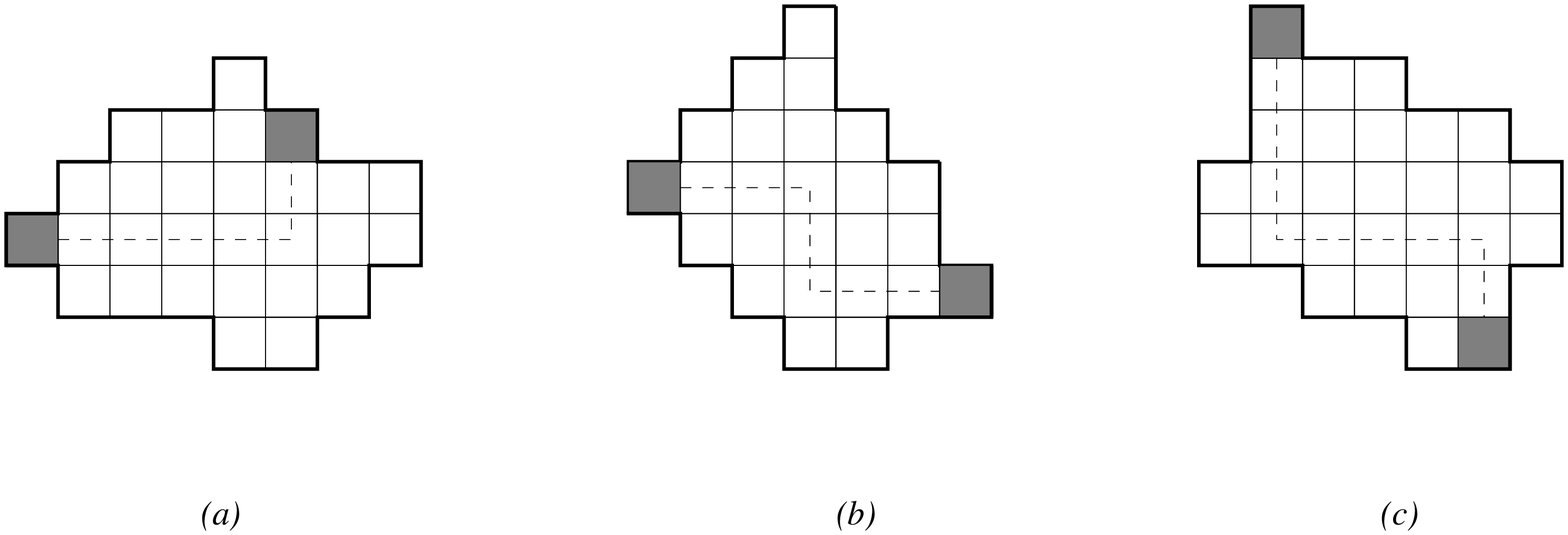}}}
\caption{(a) a \textsf{L}-convex polyomino, and a monotone path with a single
change of direction joining two of its cells; (b) a \textsf{Z}-convex but not
\textsf{L}-convex polyomino: the two highlighted cells cannot be connected by a
path with only one change of direction; (c) a centered polyomino (not
\textsf{L}-convex).}\label{lconv}
\end{figure}

The authors of~\cite{lconv} further proposed a classification of
convex polyominoes based on the number of changes of direction in the
paths connecting any two cells of a polyomino. More precisely, a
convex polyomino is $k$-convex if every pair of its cells can be
connected by a monotone path with at most $k$ changes of direction. In
a convex polyomino of the first level of this classification, any two
cells can be connected by a path with at most one change of direction:
in view of the \textsf{L}-shape of these paths, 1-convex polyominoes
are also called \textsf{L}-convex. The reader can easily check that in
Fig.~\ref{lconv}, the polyomino (a) is \textsf{L}-convex, while the
polyominoes (b), (c) are not, but are 2-convex.

This class of polyominoes has been considered from several points of
view: in~\cite{lconv2} it is shown that the set of \textsf{L}-convex
polyominoes is well-ordered with respect to the sub-picture order,
in~\cite{lconv3} the authors have investigated some tomographical
aspects of this family, and have shown that \textsf{L}-convex
polyominoes are uniquely determined by their horizontal and vertical
projections. Finally, in~\cite{enum} it is proved that the number
$g_n$ of \textsf{L}-convex polyominoes with semi-perimeter $n+2$
satisfies the recurrence relation:
\[
g_n= 4g_{n-1}-2g_{n-2}, \qquad n\geq 3,
\]
with $g_0=1$, $g_1=2$, $g_2=7$. In other terms the generating function
of \textsf{L}-convex polyominoes is rational:
\begin{eqnarray*}
G(t)&=&\sum_{n\geq0}g_nt^{n+2}=t^2+2t^3+7t^4+24t^5+82t^6+280t^7+O(t^8)\\
&=&\frac{1-2t+t^2}{1-4t+2t^2}.
\end{eqnarray*}
Indeed, in
\cite{fpsac}, the authors have provided an encoding of
\textsf{L}-convex polyominoes by words of a regular language, and have
furthermore studied the problem of enumerating \textsf{L}-convex
polyominoes with respect to the area.

In view of the definition of $\textsf{L}$-convex polyominoes as
$1$-convex polyominoes, it is natural to investigate which of the
previous properties remain true for some classes of $k$-convex
polyominoes, with $k>1$. Concerning enumeration in particular, one
would like to know if the generating functions of $k$-convex
polyominoes are rational, algebraic, or holonomic.

\subsection{\textsf{Z}-convex polyominoes}
In the present paper we deal with the family of {\em 2-convex
polyominoes}, which we rename \textsf{Z}-convex polyominoes in analogy
with the \textsf{L}-convex notation.  We shall prove the following
results:
\begin{theorem}
The generating function $P(t)$ of \textsf{Z}-convex polyominoes with
respect to the semi-perimeter is
\begin{eqnarray*}
P(t)&=&\sum_{n\geq0}p_nt^{n+2}\;=\;t^2+2t^3+7t^4+28t^5+116t^6+484t^8+O(t^8),\\
&=&
\frac{2t^4(1-2t)^2d(t)}
{(1-4t)^2(1-3t)(1-t)}+\frac{t^2(1-6t+10t^2-2t^3-t^4)}
{(1-4t)(1-3t)(1-t)},
\end{eqnarray*}
where
\[
d(t)=\frac12(1-2t-\sqrt{1-4t}).
\]
More generally, the generating function $P(x,y)$ of \textsf{Z}-convex
polyominoes with respect to the numbers of rows and columns is a
rational power series $R(x,y,d(x,y))$ in $x$, $y$ and the unique
power series $d(x,y)$ solution of the equation $d=(x+d)(y+d)$,
\[
d(x,y)=\frac12(1-x-y-\sqrt{\Delta}),
\]
where
\[
\Delta=(1-x-y)^2-4xy=(1-x-y)^2(1-\frac{4xy}{(1-x-y)^2}).
\]
More precisely,
\begin{eqnarray*}
P(x,y)&=&
\frac{2x^2y^2d(x,y)}{\Delta^2}
\frac{(1-x-y)^2}{((1-x-y)^2-xy)}
\\&&
+\;\frac{xy(1-x-y)^2(1-x-y-xy)-x^2y^2(1-x-y-3xy)}{\Delta\;((1-x-y)^2-xy)}.
\end{eqnarray*}
\end{theorem}

As conjectured by Marc Noy \cite{noy}, the asymptotic number of
\textsf{Z}-convex polyominoes with semi-perimeter $n+2$ grows like
$n\cdot 4^n$ (more precisely, $p_n\sim\frac n{24}\cdot 4^n$, so that
$f_n/p_n\to3$), while the number of \textsf{L}-convex polyominoes
grows only like $(2+\sqrt2)^n$, and the number of centered polyominoes
(see below) grows like $4^n$.

The fact that the generating function ends up in the same algebraic
extension as convex polyominoes looks surprising to us because we were
unable to derive it using the standard approaches to convex polyomino
enumeration (Temperley-like methods, wasp-waist decompositions, or
inclusion/exclusion on walks).  Instead, one interesting feature of
our paper is a construction of polyominoes by ``inflating'' smaller
one along a hook. We believe that this approach could in principle
allow for the enumeration of $k$-convex polyominoes in general.

The rest of the paper is organized as follows. The general strategy of
decomposition is explained in Section~\ref{sec:classification}. The
different cases are listed and the corresponding relations for
generating functions are derived in Section~\ref{sec:decomposition}.
Finally the resulting system of equations is solved in
Section~\ref{sec:resolution}.

\section{Classification and general strategy}\label{sec:classification}

In order to present our strategy for the decomposition, we need to
distinguish between several types of \textsf{Z}-convex polyominoes.

\subsection{Centered polyominoes} 
The first class we consider is the set $\mathcal{C}$ of {\em horizontally
centered} (or simply {\em centered}) convex polyominoes. A convex
polyomino is said to be centered if it contains at least one row
touching both the left and the right side of its minimal bounding
rectangle (see Fig.~\ref{lconv} (c)). Observe that centered
polyominoes have a simple characterization in terms of monotone paths:
\begin{lemma}\label{centered}
A convex polyomino is centered if and only if any pair of its cells can be connected by means of a path
$S^{h_1}E^{k}S^{h_2}$ or $S^{h_1}W^{k}S^{h_2}$, with $h_1,h_2,k \geq 0$.
\end{lemma}
In particular any \textsf{L}-convex polyomino is centered, and, more
importantly for us, any centered polyomino is \textsf{Z}-convex, while the
converse statements do not hold. Figure \ref{lconv} (c) shows a
centered polyomino which is not \textsf{L}-convex, and Figure \ref{lconv} (b) a
\textsf{Z}-convex polyomino which is not centered.

Centered convex polyominoes can also be described as made of two stack
polyominoes glued together at their basis. As we shall see in
Section~\ref{sec:centered}, this decomposition allows to compute
easily their generating function. In particular we shall not need here
to deal with the monotone paths.

\subsection{Non centered polyominoes} 

\begin{figure}
\centerline{\includegraphics[scale=.5]{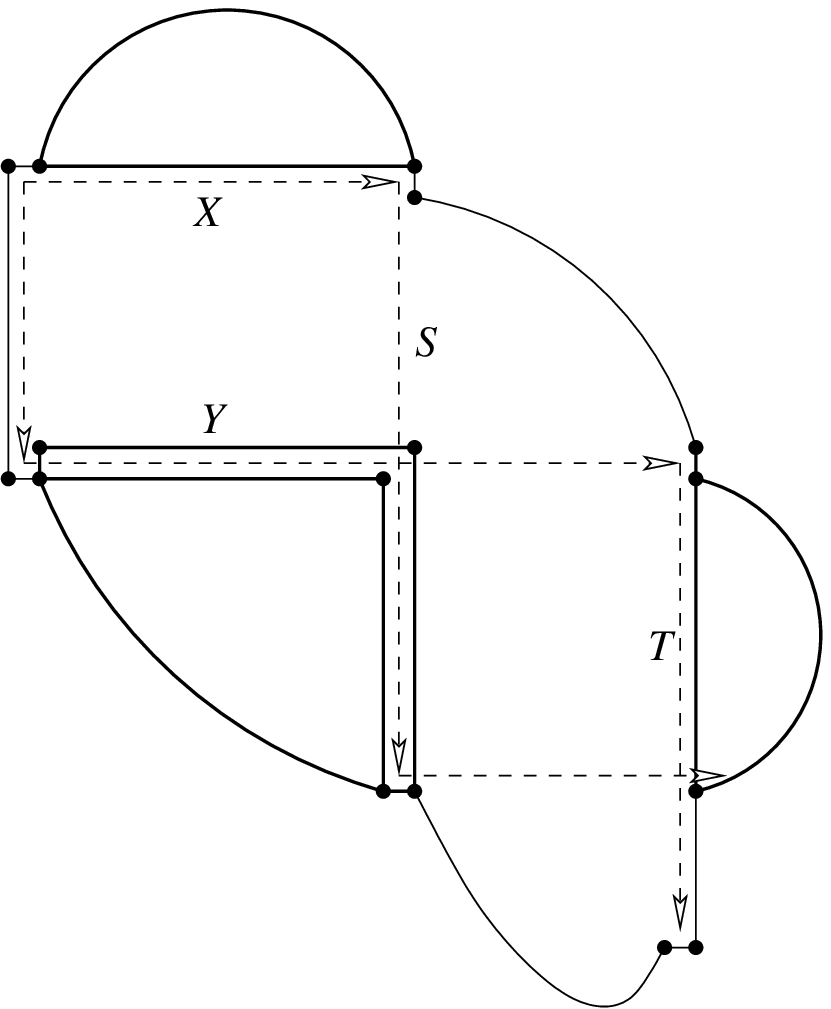}
\hspace{1em}\includegraphics[scale=.5]{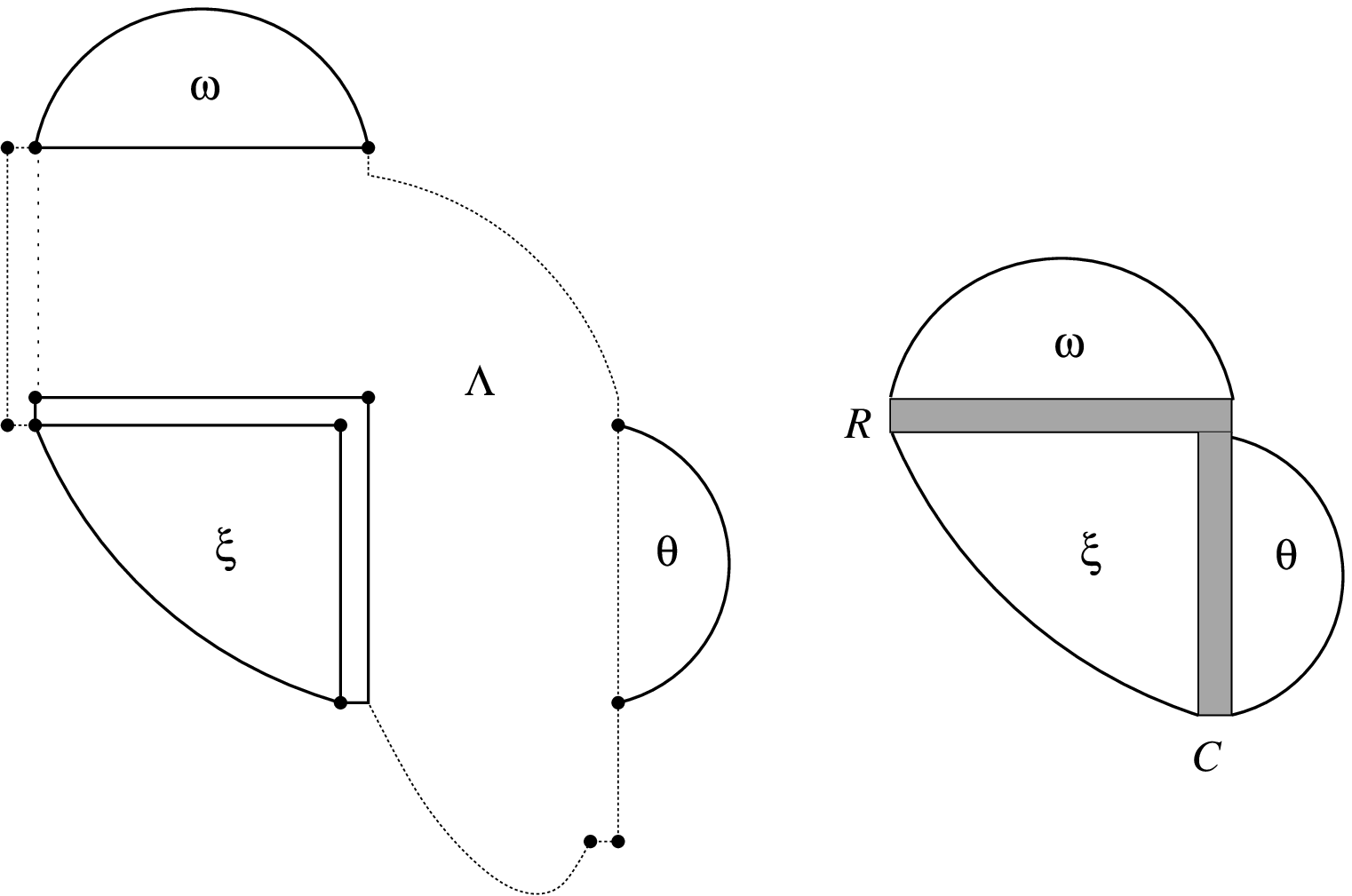}} \caption{A
non-centered polyomino of class $\mathcal D$ with its rows $x$, $y$,
and columns $s$ and $t$; its division into Regions $\omega$, $\xi$, $\theta$ and $\Lambda$; its reduction.}\label{classe1}
\end{figure}

Let us thus turn to non centered polyominoes. The starting point of
our decomposition is that we wish to remove the leftmost column. By
definition of \textsf{Z}-convexity, any two cells must be connected by a
path of type $S^{h_1}E^{k}S^{h_2}$, $S^{h_1}W^{k}S^{h_2}$,
$E^{h_1}N^{k}E^{h_2}$ or $E^{h_1}S^{k}E^{h_2}$, with $h_1,h_2,k \geq
0$. In particular, we are interested in the set of cells that can be
reached in this way from the cells of the leftmost column.

Let $P$ be a non-centered convex polyomino, and let $c_1(P)$
(briefly, $c_1$) denote its leftmost column, and let us consider the
following rows (as sketched in Fig.~\ref{classe1}):
\begin{itemize}
\item The row $X$ which contains the top cell of $c_1$.
\item The row $Y$ which contains the bottom cell of $c_1$.
\end{itemize}
Since the polyomino $P$ is convex and non-centered, its rightmost
column does not intersect any row between $X$ and $Y$, hence it is
placed entirely above $X$ or below  $Y$. 

This remark leads to the following definitions:
\begin{itemize}
\item A non-centered convex polyomino is \emph{ascending} if its rightmost
 column is above the row $X$.  Let ${\mathcal U}$ denote the set of
 descending \textsf{Z}-convex polyominoes.
\item A non-centered convex polyomino is \emph{descending} if its rightmost
column is below the row $Y$.  Let ${\mathcal D}$ denote the set of
ascending \textsf{Z}-convex polyominoes.
\end{itemize}
The whole set of \textsf{Z}-convex polyominoes is given by the union
of the three disjoint sets $\mathcal C$, $\mathcal D$, and $\mathcal
U$. Moreover, by symmetry, for any fixed size, $\mathcal D$ and
$\mathcal U$ have the same number of elements, thus, we will only
consider non-centered polyominoes of the class $\mathcal D$, as the
one represented in Fig.~\ref{classe1}.

A first property of polyominoes of class ${\mathcal D}$ is the
following consequence of their convexity: the boundary path from the
end of row $X$ to the end of row $Y$ is made only of south and east
steps.  

\subsection{The strategy}

Let us denote by $S$ and $T$ the columns starting from the rightmost
cell of $X$ and $Y$ respectively, and running until they reach the
bottom of the polyomino (see Fig.~\ref{classe1}).  The rows and
columns $X$, $Y$, $S$ and $T$ allow us to individuate four connected
sets of cells in a non-centered convex polyomino, as illustrated
by Figure~\ref{classe1}:
\begin{enumerate}
\item the set of cells strictly above $X$, called $\omega$;
\item the set of cells strictly on the right of $T$, called $\theta$;
\item the set of cells that are at the same time below $Y$ and on the
  left of $S$, called $\xi$;
\item the remaining set of cells, called $\Lambda$: these cells are
  either between $X$ and $Y$, or between $S$ and $T$ (or both).
\item[$\bullet$] In the previous definitions, the hook $H$ starting
horizontally with the left hand part of $Y$ and continuing down with
the bottom part of $S$ is included in $\xi$. The other cells of the row
$X$, $Y$ and columns $S$ and $T$ are included in $\Lambda$.
\end{enumerate}

The cells of $\theta$ require at least two turns to be reached with a
monotone path from the cells of $c_1$. The \textsf{Z}-convexity thus
induces a restriction on the position of the lowest cells of $\theta$.
\begin{property}\label{rem:lines}
The region $\theta$ of a non-centered \textsf{Z}-convex polyomino
contains no cell lower than the lowest cell of its column $S$.
\end{property}
If a row between $X$ and $Y$ reaches the right side of the bounding
box, the polyominio is centered:
\begin{property}
The set $\theta$ of a non-centered convex polyomino is non empty. 
\end{property}

As already mentioned, we wish to decompose polyominoes of $\mathcal D$
by removing the leftmost column. For the decomposition to be bijective
we then need to be able to replace a column to the left of a
polyomino. But, as the reader can verify, if one takes a
\textsf{Z}-convex polyomino and add a leftmost column, it is not so
easy to grant \emph{a priori} that Property~\ref{rem:lines} will be
satisfied by the rows and columns $X$, $Y$, $S$, and $T$ of the grown
polyomino.

In order to circumvent this problem, our decomposition will consist
into removing the whole region $\Lambda$ together with the leftmost
column.  More precisely, given a descending polyomino $P$, let us
define its \emph{reduction} $\Phi(P)$ as the polyomino obtained as
follows (see Figure~\ref{classe1}):
\begin{itemize}
\item glue region $\omega$ to $\xi$, keeping the relative abscissa of cells
between $\omega$ and $\xi$;
\item glue region $\theta$ to $\omega\cup\xi$ by keeping the relative
ordinates of cells between $\xi$ and $\theta$.
\end{itemize}
Since the hook $H$ is kept in $\Phi(P)$, $\omega$ and $\xi$ have at
least one common column (as soon as $\omega$ is non empty) and $\xi$
and $\theta$ have at least one common row, so that the reduction makes
sense and it is a polyomino, in which we highlight the hook $H$. (The
hook is highlighted in order to make easier the forthcoming description of
the inverse construction.)
 
The following lemma explains our interest in this reduction.
\begin{lemma}\label{lemma}
A descending convex polyomino is \textsf{Z}-convex if and only if it
satisfies Property~\ref{rem:lines} and its reduction $\Phi(P)$ is
\textsf{Z}-convex.
\end{lemma}
\begin{proof}
Assume first that $P$ is \textsf{Z}-convex. Then
Property~\ref{rem:lines} is satisfied and a monotone path connecting a
cell $x$ to a cell $y$ of $\Phi(P)$ can easily be constructed from the
monotone path connecting $x$ and $y$ in $P$: any section of the path
in the deleted region $\Lambda$ can be replaced by a simpler section
in the hook. 

Conversely assume that $\Phi(P)$ is \textsf{Z}-convex, that $P$
satisfies Property~\ref{rem:lines} is satisfied, and let $(x,y)$ be
two cells of $P$. If $x$ and $y$ are not in $\Lambda$ then there
exists a monotone path in $\Phi(P)$ connecting these points, and there
is no need to add a turn to extend this path into a monotone path in
$P$. If $x$ belongs to $\Lambda$, one easily construct the path in
each case $y\in\omega$, $y\in \xi$ and, using Lemma~\ref{rem:lines},
$u\in \theta$.
\end{proof}

\begin{figure}
\centerline{\hbox{\includegraphics{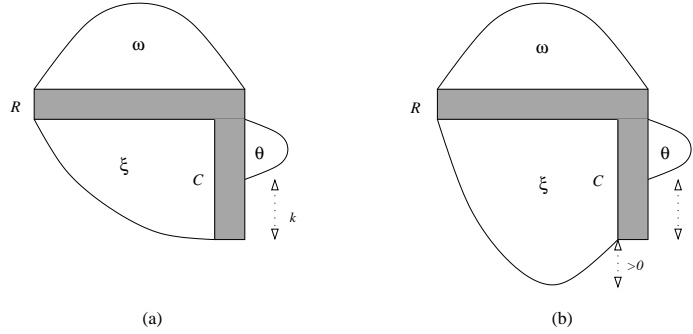}}} \caption{(a) a
hooked polyomino with hook of type $A$ and (b) one with hook of type
$B$.}\label{fig:hooked}
\end{figure}

To characterize the set of polyominoes that can occur in the image of
$\mathcal{D}$ by $\Phi$, let us define a \emph{hooked polyomino} as a
polyomino $P$ of $\mathcal{C}\cup \mathcal{D}$ in which a hook is
highlighted, in such a way that
\begin{itemize}
\item the hook is made of a top row (the \emph{arm} of hook) starting
in the leftmost column of $P$ and traversing the polyomino, and a
partial column (the \emph{leg} of the hook) starting in the right most
cell of the top row and including all cells below in this column,
\item the region on the right hand side of the hook is non empty.
\end{itemize}
The hook is called a \emph{hook of type $A$} if its bottom cell
belongs to the lowest row of $P$, and a \emph{hook of type $B$}
otherwise (see Figure~\ref{fig:hooked}). The following lemma is an
immediate consequence of the definition of hooked polyominoes.
\begin{property}
The reduction of a polyomino of $\mathcal{D}$ is a hooked polyomino.
\end{property}

In view of Lemma~\ref{lemma}, our strategy will consist in the
description of the types of region $\Lambda$ that can be added to a
hooked polyomino so that the ``inflated'' polyomino satisfies
Property~\ref{rem:lines}.

\subsection{Generating functions}
We shall compute the generating function $P(x,y)$ of \textsf{Z}-convex
polyominoes with respect to the number of columns, or \emph{width}
(variable $x$) and to the number of rows, or \emph{height} (variable
$y$). In order to do that we shall need generating functions of hooked
polyominoes with respect to the height and the width, but also with respect to
an auxiliary parameter $k$ which will be marked by a variable $u$:
given a hooked polyomino, the parameter $k$ is a non negative integer
indicating the difference of ordinate between the lowest cell of the
leg of the hook and the lowest cell of the next column to the right
(see Figure~\ref{fig:hooked}). This definition makes sense since the
region on the right hand side of the hook is assumed non empty.

We shall more precisely use the generating functions
\begin{itemize}
\item $C_A(x,y,u)$ of hooked centered polyominoes with hook of type $A$, 
\item $C_B(x,y,u)$ of hooked centered polyominoes with hook of type $B$, 
\item $A(x,y,u)=\sum_k a_k(x,y)u^k$ of hooked non-centered polyominoes with hook of type $A$, 
\item $B(x,y,u)=\sum_k b_k(x,y)u^k$ of hooked non-centered polyominoes with hook of type $B$.
\end{itemize} 
Most of the time we drop the variables $x,y$ and use the shorthand
notation $A(u)=A(x,y,u)$, $a_k=a_k(x,y)$, etc.

\section{Decompositions}\label{sec:decomposition}

\begin{figure}
\centerline{\hbox{\includegraphics[scale=.5]{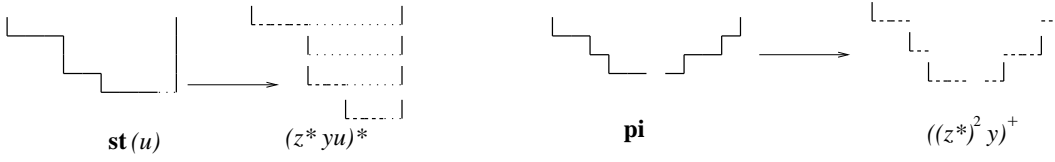}}}
\caption{The decomposition of a staircase $\mathbf{st}(u)=(z^*yu)^*$ and
of a non empty pile $\mathbf{pi}=((z^*)^2y)^+$.}\label{fig:toolbox}
\end{figure}

We shall need the following elementary notations and results,
illustrated by Figure~\ref{fig:toolbox}:
\begin{itemize}
\item The sequence notation for formal power series is
$w^*=\frac1{1-w}$. The non-empty sequence notation is $w^+=w\cdot
w^*=\frac{w}{1-w}$.
\item The generating function of possibly empty staircases with width
marked by $z$ and height marked by $yu$ is
$\mathbf{st}(u)=(z^*yu)^*$.
\item The generating function of non empty piles of lines with 
width marked by $z$ and height marked by $y$ is
$\mathbf{pi}=((z^*)^2y)^+$.
\item  
Given a generating function $F(u)=\sum_{n\geq0}f_nu^n$ we define the series
\[
F(u,v)=\sum_{n\geq0}\sum_{i+j=n}f_nu^iv^j=\frac{uF(u)}{u-v}+\frac{vF(v)}{v-u},
\]
where the inner summation is on non negative $i$ and $j$ with $i+j=n$, and 
\[
F(u,v,w)=\sum_{n\geq0}\sum_{i+j+k=n}f_nu^iv^jw^k=\frac{u^2F(u)}{(u-v)(u-w)}+\frac{v^2F(v)}{(v-u)(w-u)}+\frac{w^2F(w)}{(w-u)(w-v)},
\]
where the inner summation is on non negative $i$, $j$ and $k$ with
$i+j+k=n$.  The values at $u=v$ of the previous series can be obtained
by continuity:
\[
A(u,u,w)=\left(\frac{u^2A(u)}{u-w}\right)'+\frac{w^2A(w)}{(w-u)^2}=
\frac{u^2A'(u)}{u-w}+\frac{u(u-2w)A(u)}{(u-w)^2}+\frac{w^2A(w)}{(w-u)^2}.
\]
\end{itemize}


\subsection{Centered polyominoes}\label{sec:centered}

Recall that a centered polyomino is a polyomino that contains at least
one row touching both the left and the right hand side of its minimal
bounding rectangle. We need to count polyominoes of the family
$\mathcal{C}$ of centered polyominoes but also of the families
$\mathcal{C}_A$ and $\mathcal{C}_B$ of hooked polyominoes with a hook
of type $A$ and $B$ respectively.

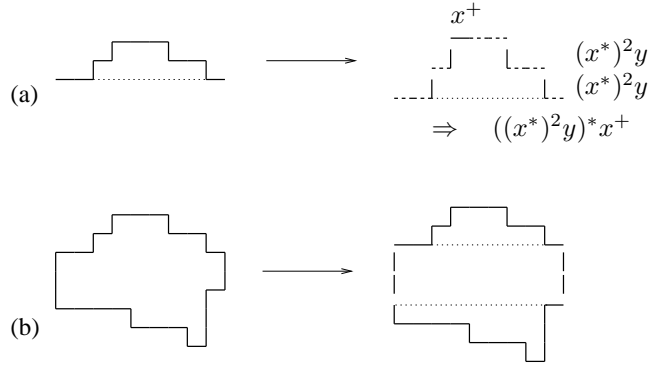
\begin{figure}
\centerline{
\input{double-stack.pstex_t}
}
\caption{(a) The decomposition of a stack polyomino with baseline
width marked by $x$ and height marked by $y$. (b) The decomposition of
a centered polyomino into a non-empty sequence of central rows and two
stack polyominoes.}
\label{fig:double-stack}
\end{figure}

Let $S(x,y)$ be the generating function of stack polyominoes with $x$
marking the length of the baseline and $y$ marking the height.
In view of Figure~\ref{fig:double-stack}(a),
\[
S(x,y)=x^+\cdot((z^*)^2y)^*=\frac{x(1-x)}{(1-x)^2-y}
=\left(\frac12\frac1{1-\frac x{1-\sqrt y}}+
\frac12\frac1{1-\frac x{1+\sqrt y}}-1\right).
\]
Observe then that for any power series $F(x,y)$ the Hadamard product
$S(x,y)\odot_x F(x,y)$ is equal to:
\[
\frac12 \left(F\big(\frac{x}{1-\sqrt y},y\big)+F\big(\frac{x}{1+\sqrt y},y\big)\right)-F(0),
\]
which is a rational function of $x$ and $y$ if $F(x)$ is.

In view of Figure~\ref{fig:double-stack}(b), centered polyominoes are
formed of a centered rectangle supporting 2 strictly smaller stacks
polyominoes:
\[
C(x,y)=y^+[S^>(x,y)\odot_xS^>(x,y)],
\]
where $S^>(x,y)$ stands for the generating function of stack
polyominoes with a first row strictly smaller than the baseline (so
that the central rectangle is effectively given by the factor $y^+$).
The series $S^>(x,y)$ is readily obtain by difference,
\[
S^{>}(x,y)=S(x,y)-yS(x,y)=\frac{x(1-x)(1-y)}{(1-x)^2-y},
\]
and computing the Hadamard product with the previous formula yields:
\[
C(x,y)=\frac{xy(-y-xy+1-2x+x^2)(1-y)}{(1-x-y)(x^2-2xy-2x+y^2-2y+1)}.
\]

\begin{figure}
\input{classeCA.pstex_t}
\caption{Construction of elements of the classes $C_A(u)$ and $C_B(u)$.}
\label{fig:centered-hook}
\end{figure}
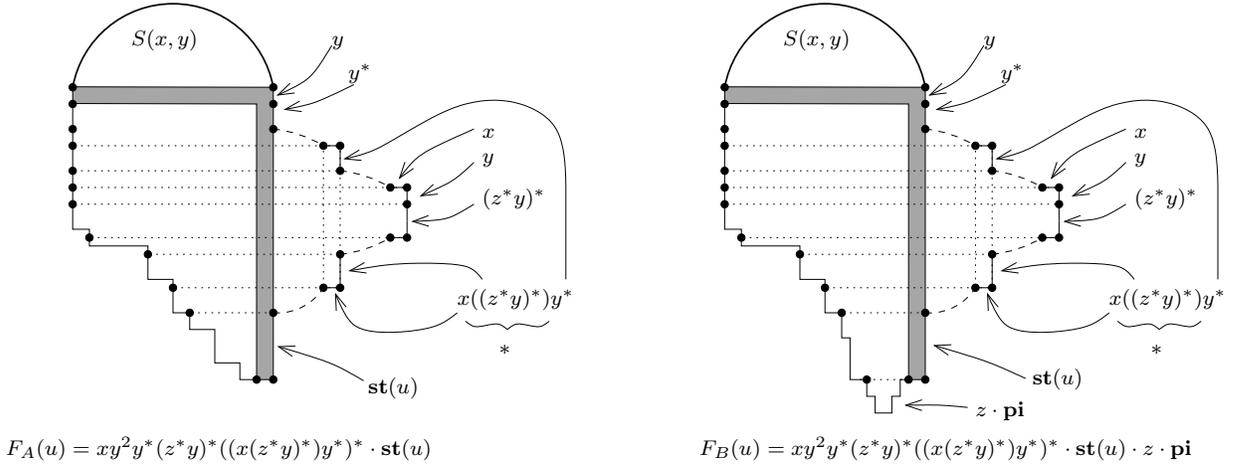

We shall also need centered polyominoes with a marked hook, of type
$A$ and $B$.  As illustrated by Figure~\ref{fig:centered-hook}, the
series for the first type is $C_A(u)=S(z,y)\odot_z F_A(zx,x,y,u)$
where
\[F_A(z,x,y,u)=y^2z^+\cdot\mathbf{st}(u)\cdot(xy^*(yz^*)^*)^+.\]
Indeed, with $z$ marking columns on the left hand side of the leg of
the hook, the Hadamard product accounts for gluying, along the arm of
the hook, a staircase, with generating series $S(z,y)$, to the rest of
the polyomino, with generating series $F_A(z,x,y,u)$: in this later
series, a factor $xy(z^*y)^*$ corresponds to the central rectangle;
each factor $x((z^*y)^*)y^*$ corresponds to a column on the right of
the hook and to the lines having their rightmost cell in that column;
the factor $\mathbf{st}(u)$ corresponds to the bottom staircase made
of lines having their rightmost cell in the hook.

Similarly, the series for the second type is $C_B(u)=S(z,y)\odot_z
F_B(zx,x,y,u)$ where
\[F_B(z,x,y,u)=F_A(z,x,y,u)\cdot z\cdot\mathbf{pi},\]
with the extra factor corresponding to cells lower than the leg of the
hook.

\subsection{Hooked polyominoes with hooks of type $A$}

\begin{figure}
\begin{center}
\input{classeA1-3.pstex_t}
\end{center}
\caption{Construction of elements of classes $A_1$, $A_2$ and $A_3$,
from the hooked polyomino of
Figure~\ref{fig:hooked}(a).}\label{classeA1-3}

\bigskip
\begin{center}
\input{classeA4-5.pstex_t}
\end{center}
\caption{Construction of elements of classes $A_4$, and
$A_5$, from the hooked polyomino of
Figure~\ref{fig:hooked}(a).}\label{classeA4-5}
\end{figure}

A hooked polyomino with hook ok type $A$ can be a hooked centered
polyomino (with gf $C_A(u)$ already computed) or can be obtained from
its reduction which must be a hooked polyomino with a hook of type $A$
(recall that type $A$ means that the leg of the hook reaches the lowest
row of the polyomino). Let us describe the different cases, with
respect to the properties of the resulting inflated polyomino:
\begin{itemize}
\item 
The leg of the hook and the two columns $S$ and $T$ have
same abscissa (Figure~\ref{classeA1-3}, left): let
\[
A_1(u)={x(y^*)^2}\cdot A(u).
\]
\item 
The leg of the hook has same abscissa as the column $T$ but not as $S$
(Figure~\ref{classeA1-3}, middle): by definition of type $A$, the
column $S$ cannot be longer than the leg of the hook, and
\[
A_2(u)={xy^*y^+}\cdot\mathbf{st}(u)\cdot z^+A(u).
\]
The series $A_2(u)$ apparently does not takes into account the
construction of the staircase starting on the righthand side of the
column $S$ and connecting it to the top-right angle of the
hook. Instead each column between column $S$ (excluded) and the leg of
the hook (included) is marked by a factor $z$. However upon setting
$z=xy^*$, each column marked by $z$ gets a factor $x$ and a factor
$y^*$ that accounts for the rows ending in that column. The generating function
of polyominoes of this case is thus $A_2(u)\big|_{z=xy^*}$.

In all forthcoming cases, we describe similarly generating functions
of polyominoes without the staircase connecting $X$ to the top-right
corner of the hook. In other terms, in the following pages $z$ is to be understood as a shorthand notation for $xy^*$.

\item 
The leg of the hook has same abscissa as the column $S$ but not as $T$
(Figure~\ref{classeA1-3}, right): by definition of type $A$, the leg
of the hook is at least as long as the column $T$, and
\[
A_3(u)={xy^*y^+}\cdot z^+A(u,z^*).
\]
As suggested by Figure~\ref{classeA1-3}, the factor $A(u,z^*)$
accounts for the fact that from a polyomino of type $A$ with parameter
$k$, one constructs a new polyomino with parameter $i$ with $0\leq
i\leq k$.
\item 
The abscissa of the leg of the hook is strictly between $S$ and $T$
and the cells marked by a factor $u$ are strictly below the lowest cells of
columns $S$ and $T$  (Figure~\ref{classeA4-5}, left):
\[
A_4(u)={x(y^+)^2}\cdot\mathbf{st}(u)\cdot \mathbf{pi}\cdot(z^+)^2A(z^*).
\]
\item 
The abscissa of the leg of the hook is strictly between $S$ and $T$
and the cells marked by a factor $u$ intersects the baseline of column
$S$ (Figure~\ref{classeA4-5}, right):
\[
A_5(u)={x(y^+)^2}\cdot \mathbf{st}(u)\cdot(z^+)^2A(u,z^*).
\]
\end{itemize}
The generating series of hooked polyominoes with a hook of type $A$ is
then
\[
A(u)=C_A(u)+\sum_{i=1}^5A_i(u)\big|_{z=xy^*}.
\]


\subsection{Hooked polyominoes with hooks of type $B$}
\begin{figure}
\begin{center}
\mbox{\input{classeB1.pstex_t}
\hspace{1em}
\input{classeB4.pstex_t}
}
\end{center}
\caption{Construction of elements of classes $B_1$, and
$B_4$, from the hooked polyomino of
Figure~\ref{fig:hooked}(a).}\label{classeB1-4}

\bigskip
\begin{center}
\input{classeB5-7.pstex_t}
\end{center}
\caption{Construction of elements of classes $B_5$, $B_6$ and
$B_7$, from the hooked polyomino of
Figure~\ref{fig:hooked}(b).}\label{classeB5-7}
\end{figure}

A hooked polyomino with hook of type $B$ can be a hooked centered
polyomino (with gf $C_B(u)$ already computed) or can be obtained by
inflating its reduced polyomino. We start with those that are produced
from a hooked polyomino with hook of type $A$, and we give again the different
cases with respect to the properties of the obtained inflated
polyomino:
\begin{itemize}
\item 
The leg of the hook has same abscissa as column $T$ and it is strictly
longer than column $S$: in order to produce a hook of type $B$, some
other column between $S$ and $T$ must be even longer, and
\[
B_1(u)=A_2(u)\cdot z\cdot\mathbf{pi}={xy^*y^+}\cdot\mathbf{pi}\cdot z\cdot \mathbf{st}(u)\cdot z^*yu\cdot z^+A(u).
\]
In agreement with the relation $B_1(u)=A_2(u)\cdot z\cdot\mathbf{pi}$,
polyominoes of $B_1$ (Figure~\ref{classeB1-4}, left) can be obtained
from polyominoes of $A_2$ (Figure~\ref{classeA1-3}, middle) by adding
a non-empty pile of rows just before the leg of the hook.
\item 
The abscissa of the leg of the hook is strictly between $S$ and $T$,
and the cells marked by a factor $u$ are strictly below the lowest
cells of columns $S$ and $T$:
\[
B_2(u)=A_4(u)\cdot z\cdot\mathbf{pi}={x(y^+)^2}\cdot\mathbf{pi}\cdot z\cdot\mathbf{st}(u)\cdot \mathbf{pi}\cdot(z^+)^2A(z^*).
\]
These polyominoes are obtained from the polyominoes of $A_4$ 
upon adding a non-empty pile of rows just before the leg of the hook.
\item 
The abscissa of the leg of the hook is strictly between $S$ and $T$,
and the cells marked by a factor $u$ intersects the baseline of column
$S$.
\[
B_3(u)=A_5(u)\cdot z\cdot\mathbf{pi}={x(y^+)^2}\cdot\mathbf{pi}\cdot z\cdot \mathbf{st}(u)\cdot(z^+)^2A(u,z^*).
\]
These polyominoes are obtained from the polyominoes of $A_5$ upon
adding a non-empty pile of rows just before the leg of the hook.
\item 
The abscissa of the leg of the hook is strictly between $S$ and $T$,
and the cells marked by a factor $u$ are strictly above the lowest
cell of column $S$.
\[
B_4(u)={x(y^+)^2}\cdot(1+\mathbf{pi}\cdot z)(z^+)^2(A(z^*,z^*,u)-A(z^*,u)).
\]
Observe that difference is due to the restriction $j\neq0$, as
illustrated by the Figure~\ref{classeB1-4}: the leg of the hook must end
strictly above the lowest cell of columns $S$, so that one must have
$j\geq1$.
\end{itemize}
Now we present the cases produced from a hooked polyomino with hook of
type $B$, again arranged according to the properties of the resulting
polyomino. Observe that in these cases the column $S$ is at least as
long as the column $T$:
\begin{itemize}
\item 
The leg of the hook and the columns $S$ and $T$ have the same
abscissa:
\[
B_5(u)=x(y^*)^2\cdot B(u)
\]
\item 
The leg of the hook has the same abscissa as the
column $S$ or the same abscissa as the column $T$ (but not both):
\[
B_6(u)=2\cdot{xy^*y^+}\cdot z^*B(z^*,u).
\]
\item 
The abscissa of the leg of the hook is strictly between column $S$ and $T$:
\[
B_7(u)={x(y^+)^2}\cdot(z^+)^2B(z^*,z^*,u).
\]
\end{itemize}
The generating series of hooked polyominoes with a hook of type $B$ is
then
\[
B(u)=C_B(u)+\sum_{i=1}^7B_i(u)\big|_{z=xy^*}.
\]

\subsection{\textsf{Z}-convex polyominoes}

\begin{figure}
\begin{center}
\input{classeP1-3.pstex_t}
\end{center}
\caption{Construction of elements of classes $P_1$, $P_2$ and $P_3$,
from the hooked polyomino of
Figure~\ref{fig:hooked}(a).}\label{classeP1-3}

\bigskip
\begin{center}
\input{classeP4-5.pstex_t}
\end{center}
\caption{Construction of elements of classes $P_4$, and
$P_5$, from the hooked polyomino of
Figure~\ref{fig:hooked}(b).}\label{classeP4-5}
\end{figure}

Again we start with polyominoes that are produced from a hooked polyomino with hook of type $A$:
\begin{itemize}
\item 
The columns $S$ and $T$ have the same abscissa:
\[
P_1=xy^*A(1).
\]
\item 
The columns $S$ and $T$ have distinct abscissa and the column $S$ is
strictly shorter than $T$:
\[
P_2={xy^+}\cdot\mathbf{pi}\cdot zA(1).
\]
\item 
The columns $S$ and $T$ have distinct abscissa and the
column $S$ is at least as long as $T$:
\[
P_3={xy^+}\cdot(1+z\cdot\mathbf{pi})(z^*A(z^*)-A(1)).
\]
The difference is due to the fact that at least one horizontal column
must be inserted at the level of the rows that were marked by the
factor $u$ to ensure that the column $T$ is not longer that $S$.
\end{itemize}
Next we present the polyominoes  obtained from a hooked
polyomino with hook of type $B$:
\begin{itemize}
\item 
The columns $S$ and $T$ have the same abscissa:
\[
P_4=xy^*B(1).
\]
\item 
 The columns $S$ and $T$ have distinct abscissa and the
lowest cell of $S$ is below or at the same level as the lowest cell of
$T$:
\[
P_5=xy^+\cdot(z^*B(z^*)-B(1)).
\]
\end{itemize}
Finally the generating function of \textsf{Z}-convex polyominoes is
\[
P=C+\sum_{i=1}^5P_i\big|_{z=xy^*}.
\]

\section{Resolution}\label{sec:resolution}
In view of the previous section, upon setting as announced $z=xy^*$,
the system of equations defining the series $P$ has the following
form:
\begin{eqnarray*}
A(u)&=&C_A(u)+a_1(u)A(u)+a_2(u)A(z^*)\\
B(u)&=&C_B(u)+b_1(u)A(u)+b_2(u)A(z^*)+b_3(u)A'_u(z^*)+b_4(u)B(u)
+b_5(u)B(z^*)+b_6(u)B'_u(z^*)\\
P&=&C_0+p_1A(z^*)+p_2A(1)+p_3B(z^*)+p_4B(1),
\end{eqnarray*}
where the $C_A(u)$, $C_B(u)$ and $C_0$ are the rational generating
series of centered polyominoes computed in Section~\ref{sec:centered},
the $a_i(u)$ and $b_i(u)$ are explicit rational functions of $x$, $y$
and $u$, and the $p_i$ are explicit rational functions of $x$ and $y$.

The first step of the resolution is to apply the kernel method to the
first equation, which involves only $A(u)$ and $A(z^*)$ as unknown.
The kernel equation $1-a_1(u)=0$ contains a factor that can be written
\[
u=1+(x-y)u+yu^2,
\]
so that it clearly admits a power series root $c(x,y)$,
which is a refinement of the Catalan generating function 
\[
c(t,t)=\frac{1-\sqrt{1-4t}}{2t}.
\] 
Setting $u=c(x,y)$ in the first equation, the kernel is canceled and
$A(z^*)$ is obtained as
\[
A(z^*)=\frac{C_A(c)}{a_2(c)}.
\]
Then, using again the first equation of the system we derive
$A(u)$. Once $A(u)$ is known, $A'_u(z^*)$ can also be computed.

The second step consists in applying now the kernel method to the
second equation of the system, which now has three unknowns $B(u)$,
$B(z^*)$ and $B'_u(z^*)$. The kernel $1-b_4(u)$ admits two roots $R_1$
and $R_2$ that are rational power series in $x^{1/2}$ and
$y$:
\[
R_1=\frac{1-y+x^{1/2}}{1-y-x+(1-x)x^{1/2}},
\;\textrm{ and }\; 
R_2=\frac{1-y-x^{1/2}}{1-y-x-(1-x)x^{1/2}}.
\]
Using these two roots we write two linear equations for $B(z^*)$ and
$B'_u(z^*)$ and solve the system.  The resulting series are rational
series in $x$, $y$ and $C(x,y)$ (in particular fractional powers of
$x$ cancel, as one could expect from the symmetry with respect to $\pm
x^{1/2}$). Returning to the second equation of the system, we obtain
$B(u)$ and finally, turning to the third equation, the generating
function $P$ of \textsf{Z}-convex polyominoes.

It should be remarked that our method leads to heavy computations in
the intermediary steps, involving big rational expressions. The fact
that things dramatically simplify when all pieces are put together in
$P$ calls for a simpler, more combinatorial, derivation. In
particular, the expression are nicer in terms of the more symmetric
parametrization $d(x,y)=y(c(x,y)-1)$ satisfying $d=(x+d)(y+d)$.

\subsection*{Acknowledgments.} 
The authors wish to thank Andrea Frosini and Marc Noy for inspiring
discussions on the topic of this paper. Andrea Frosini also provided
the first terms of the series by exhaustive generation, allowing us to
double check our results.

\end{document}

%% file: double-stack.pstex_t
\begin{picture}(0,0)%
\includegraphics{double-stack.pstex}%
\end{picture}%
\setlength{\unitlength}{2072sp}%
\begingroup\makeatletter\ifx\SetFigFont\undefined%
\gdef\SetFigFont#1#2#3#4#5{%
  \reset@font\fontsize{#1}{#2pt}%
  \fontfamily{#3}\fontseries{#4}\fontshape{#5}%
  \selectfont}%
\fi\endgroup%
\begin{picture}(7980,4264)(9811,-9533)
\put(15076,-5461){\makebox(0,0)[lb]{\smash{{\SetFigFont{10}{12.0}{\rmdefault}{\mddefault}{\itdefault}{\color[rgb]{0,0,0}$x^+$}%
}}}}
\put(16561,-5911){\makebox(0,0)[lb]{\smash{{\SetFigFont{10}{12.0}{\rmdefault}{\mddefault}{\itdefault}{\color[rgb]{0,0,0}$(x^*)^2y$}%
}}}}
\put(14851,-6811){\makebox(0,0)[lb]{\smash{{\SetFigFont{10}{12.0}{\rmdefault}{\mddefault}{\itdefault}{\color[rgb]{0,0,0}$\Rightarrow\quad((x^*)^2y)^*x^+$}%
}}}}
\put(16561,-6316){\makebox(0,0)[lb]{\smash{{\SetFigFont{10}{12.0}{\rmdefault}{\mddefault}{\itdefault}{\color[rgb]{0,0,0}$(x^*)^2y$}%
}}}}
\end{picture}%

%% file: classeCA.pstex_t
\begin{picture}(0,0)%
\includegraphics{classeCA.pstex}%
\end{picture}%
\setlength{\unitlength}{2763sp}%
\begingroup\makeatletter\ifx\SetFigFont\undefined%
\gdef\SetFigFont#1#2#3#4#5{%
  \reset@font\fontsize{#1}{#2pt}%
  \fontfamily{#3}\fontseries{#4}\fontshape{#5}%
  \selectfont}%
\fi\endgroup%
\begin{picture}(11499,4124)(7501,-6644)
\put(16276,-2911){\makebox(0,0)[lb]{\smash{{\SetFigFont{8}{9.6}{\rmdefault}{\mddefault}{\updefault}{\color[rgb]{0,0,0}$y$}%
}}}}
\put(16426,-3211){\makebox(0,0)[lb]{\smash{{\SetFigFont{8}{9.6}{\rmdefault}{\mddefault}{\updefault}{\color[rgb]{0,0,0}$y^{*}$}%
}}}}
\put(17626,-3961){\makebox(0,0)[lb]{\smash{{\SetFigFont{8}{9.6}{\rmdefault}{\mddefault}{\updefault}{\color[rgb]{0,0,0}$y$}%
}}}}
\put(17626,-4336){\makebox(0,0)[lb]{\smash{{\SetFigFont{8}{9.6}{\rmdefault}{\mddefault}{\updefault}{\color[rgb]{0,0,0}$(z^*y)^*$}%
}}}}
\put(17626,-3736){\makebox(0,0)[lb]{\smash{{\SetFigFont{8}{9.6}{\rmdefault}{\mddefault}{\updefault}{\color[rgb]{0,0,0}$x$}%
}}}}
\put(14476,-2911){\makebox(0,0)[lb]{\smash{{\SetFigFont{8}{9.6}{\rmdefault}{\mddefault}{\updefault}{\color[rgb]{0,0,0}$S(x,y)$}%
}}}}
\put(17776,-5686){\makebox(0,0)[lb]{\smash{{\SetFigFont{8}{9.6}{\rmdefault}{\mddefault}{\updefault}{\color[rgb]{0,0,0}$*$}%
}}}}
\put(17401,-5236){\makebox(0,0)[lb]{\smash{{\SetFigFont{8}{9.6}{\rmdefault}{\mddefault}{\updefault}{\color[rgb]{0,0,0}$x((z^*y)^{*})y^*$}%
}}}}
\put(7501,-6586){\makebox(0,0)[lb]{\smash{{\SetFigFont{8}{9.6}{\rmdefault}{\mddefault}{\updefault}{\color[rgb]{0,0,0}$F_A(u)=xy^{2}y^{*}(z^{*}y)^{*}((x(z^*y)^*)y^*)^*\cdot\mathbf{st}(u)$}%
}}}}
\put(13726,-6586){\makebox(0,0)[lb]{\smash{{\SetFigFont{8}{9.6}{\rmdefault}{\mddefault}{\updefault}{\color[rgb]{0,0,0}$F_B(u)=xy^{2}y^{*}(z^{*}y)^{*}((x(z^*y)^*)y^*)^*\cdot\mathbf{st}(u)\cdot z\cdot\mathbf{pi}$}%
}}}}
\put(10426,-2911){\makebox(0,0)[lb]{\smash{{\SetFigFont{8}{9.6}{\rmdefault}{\mddefault}{\updefault}{\color[rgb]{0,0,0}$y$}%
}}}}
\put(10576,-3211){\makebox(0,0)[lb]{\smash{{\SetFigFont{8}{9.6}{\rmdefault}{\mddefault}{\updefault}{\color[rgb]{0,0,0}$y^{*}$}%
}}}}
\put(11776,-3961){\makebox(0,0)[lb]{\smash{{\SetFigFont{8}{9.6}{\rmdefault}{\mddefault}{\updefault}{\color[rgb]{0,0,0}$y$}%
}}}}
\put(11776,-4336){\makebox(0,0)[lb]{\smash{{\SetFigFont{8}{9.6}{\rmdefault}{\mddefault}{\updefault}{\color[rgb]{0,0,0}$(z^*y)^*$}%
}}}}
\put(11776,-3736){\makebox(0,0)[lb]{\smash{{\SetFigFont{8}{9.6}{\rmdefault}{\mddefault}{\updefault}{\color[rgb]{0,0,0}$x$}%
}}}}
\put(8626,-2911){\makebox(0,0)[lb]{\smash{{\SetFigFont{8}{9.6}{\rmdefault}{\mddefault}{\updefault}{\color[rgb]{0,0,0}$S(x,y)$}%
}}}}
\put(11926,-5686){\makebox(0,0)[lb]{\smash{{\SetFigFont{8}{9.6}{\rmdefault}{\mddefault}{\updefault}{\color[rgb]{0,0,0}$*$}%
}}}}
\put(11551,-5236){\makebox(0,0)[lb]{\smash{{\SetFigFont{8}{9.6}{\rmdefault}{\mddefault}{\updefault}{\color[rgb]{0,0,0}$x((z^*y)^{*})y^*$}%
}}}}
\put(10751,-6036){\makebox(0,0)[lb]{\smash{{\SetFigFont{8}{9.6}{\rmdefault}{\mddefault}{\updefault}{\color[rgb]{0,0,0}$\mathbf{st}(u)$}%
}}}}
\put(16201,-6241){\makebox(0,0)[lb]{\smash{{\SetFigFont{8}{9.6}{\rmdefault}{\mddefault}{\updefault}{\color[rgb]{0,0,0}$z\cdot\mathbf{pi}$}%
}}}}
\put(16711,-5956){\makebox(0,0)[lb]{\smash{{\SetFigFont{8}{9.6}{\rmdefault}{\mddefault}{\updefault}{\color[rgb]{0,0,0}$\mathbf{st}(u)$}%
}}}}
\end{picture}%

%% file: classeA1-3.pstex_t
\begin{picture}(0,0)%
\includegraphics{classeA1-3.pstex}%
\end{picture}%
\setlength{\unitlength}{2763sp}%
\begingroup\makeatletter\ifx\SetFigFont\undefined%
\gdef\SetFigFont#1#2#3#4#5{%
  \reset@font\fontsize{#1}{#2pt}%
  \fontfamily{#3}\fontseries{#4}\fontshape{#5}%
  \selectfont}%
\fi\endgroup%
\begin{picture}(11448,5174)(6226,-7469)
\put(12286,-3481){\makebox(0,0)[lb]{\smash{{\SetFigFont{8}{9.6}{\rmdefault}{\mddefault}{\updefault}{\color[rgb]{0,0,0}$z=xy^{*}$}%
}}}}
\put(11566,-2941){\makebox(0,0)[lb]{\smash{{\SetFigFont{8}{9.6}{\rmdefault}{\mddefault}{\updefault}{\color[rgb]{0,0,0}$y$}%
}}}}
\put(11656,-3256){\makebox(0,0)[lb]{\smash{{\SetFigFont{8}{9.6}{\rmdefault}{\mddefault}{\updefault}{\color[rgb]{0,0,0}$y^{*}$}%
}}}}
\put(12556,-4201){\makebox(0,0)[lb]{\smash{{\SetFigFont{8}{9.6}{\rmdefault}{\mddefault}{\updefault}{\color[rgb]{0,0,0}$y^{*}$}%
}}}}
\put(8866,-2581){\makebox(0,0)[lb]{\smash{{\SetFigFont{8}{9.6}{\rmdefault}{\mddefault}{\updefault}{\color[rgb]{0,0,0}$x$}%
}}}}
\put(8356,-4111){\makebox(0,0)[lb]{\smash{{\SetFigFont{8}{9.6}{\rmdefault}{\mddefault}{\updefault}{\color[rgb]{0,0,0}$y^{*}$}%
}}}}
\put(8356,-3436){\makebox(0,0)[lb]{\smash{{\SetFigFont{8}{9.6}{\rmdefault}{\mddefault}{\updefault}{\color[rgb]{0,0,0}$y^{*}$}%
}}}}
\put(12901,-7111){\makebox(0,0)[lb]{\smash{{\SetFigFont{8}{9.6}{\rmdefault}{\mddefault}{\updefault}{\color[rgb]{0,0,0}$A_{3}(u)=xy(y^{*})^{2}z\sum_{k}a_{k}\sum_{i}u^{k-i}(z^{*})^{i+1}$}%
}}}}
\put(13576,-7411){\makebox(0,0)[lb]{\smash{{\SetFigFont{8}{9.6}{\rmdefault}{\mddefault}{\updefault}{\color[rgb]{0,0,0}$=xy(y^{*})^{2}z^{+}A(u,z^{*})$}%
}}}}
\put(6376,-6736){\makebox(0,0)[lb]{\smash{{\SetFigFont{8}{9.6}{\rmdefault}{\mddefault}{\updefault}{\color[rgb]{0,0,0}$A_{1}(u)=x(y^{*})^{2}A(u)$}%
}}}}
\put(10201,-6436){\makebox(0,0)[lb]{\smash{{\SetFigFont{8}{9.6}{\rmdefault}{\mddefault}{\updefault}{\color[rgb]{0,0,0}$z^{*}$}%
}}}}
\put(11401,-7036){\makebox(0,0)[lb]{\smash{{\SetFigFont{8}{9.6}{\rmdefault}{\mddefault}{\updefault}{\color[rgb]{0,0,0}$z$}%
}}}}
\put(6226,-2611){\makebox(0,0)[lb]{\smash{{\SetFigFont{8}{9.6}{\rmdefault}{\mddefault}{\updefault}{\color[rgb]{0,0,0}$x$}%
}}}}
\put(12766,-2581){\makebox(0,0)[lb]{\smash{{\SetFigFont{8}{9.6}{\rmdefault}{\mddefault}{\updefault}{\color[rgb]{0,0,0}$x$}%
}}}}
\put(15601,-3931){\makebox(0,0)[lb]{\smash{{\SetFigFont{8}{9.6}{\rmdefault}{\mddefault}{\updefault}{\color[rgb]{0,0,0}$y^{*}$}%
}}}}
\put(15601,-3481){\makebox(0,0)[lb]{\smash{{\SetFigFont{8}{9.6}{\rmdefault}{\mddefault}{\updefault}{\color[rgb]{0,0,0}$y^{*}$}%
}}}}
\put(15466,-2941){\makebox(0,0)[lb]{\smash{{\SetFigFont{8}{9.6}{\rmdefault}{\mddefault}{\updefault}{\color[rgb]{0,0,0}$y$}%
}}}}
\put(16186,-4111){\makebox(0,0)[lb]{\smash{{\SetFigFont{8}{9.6}{\rmdefault}{\mddefault}{\updefault}{\color[rgb]{0,0,0}$z=xy^{*}$}%
}}}}
\put(16351,-6061){\makebox(0,0)[lb]{\smash{{\SetFigFont{8}{9.6}{\rmdefault}{\mddefault}{\updefault}{\color[rgb]{0,0,0}$(z^{*})^{i+1}$}%
}}}}
\put(15751,-6136){\makebox(0,0)[lb]{\smash{{\SetFigFont{8}{9.6}{\rmdefault}{\mddefault}{\updefault}{\color[rgb]{0,0,0}$z$}%
}}}}
\put(15676,-6436){\makebox(0,0)[lb]{\smash{{\SetFigFont{8}{9.6}{\rmdefault}{\mddefault}{\updefault}{\color[rgb]{0,0,0}$u^{k-i}$}%
}}}}
\put(8356,-5876){\makebox(0,0)[lb]{\smash{{\SetFigFont{8}{9.6}{\rmdefault}{\mddefault}{\updefault}{\color[rgb]{0,0,0}$u^{k}$}%
}}}}
\put(12531,-6336){\makebox(0,0)[lb]{\smash{{\SetFigFont{8}{9.6}{\rmdefault}{\mddefault}{\updefault}{\color[rgb]{0,0,0}$\mathbf{st}(u)$}%
}}}}
\put(12576,-5836){\makebox(0,0)[lb]{\smash{{\SetFigFont{8}{9.6}{\rmdefault}{\mddefault}{\updefault}{\color[rgb]{0,0,0}$u^{k}$}%
}}}}
\put(9346,-7311){\makebox(0,0)[lb]{\smash{{\SetFigFont{8}{9.6}{\rmdefault}{\mddefault}{\updefault}{\color[rgb]{0,0,0}$A_{2}(u)=xy(y^{*})^{2}z^{+}\mathbf{st}(u)A(u)$}%
}}}}
\end{picture}%

%% file: classeA4-5.pstex_t
\begin{picture}(0,0)%
\includegraphics{classeA4-5.pstex}%
\end{picture}%
\setlength{\unitlength}{2763sp}%
\begingroup\makeatletter\ifx\SetFigFont\undefined%
\gdef\SetFigFont#1#2#3#4#5{%
  \reset@font\fontsize{#1}{#2pt}%
  \fontfamily{#3}\fontseries{#4}\fontshape{#5}%
  \selectfont}%
\fi\endgroup%
\begin{picture}(11118,5249)(7741,-7544)
\put(17041,-6181){\makebox(0,0)[lb]{\smash{{\SetFigFont{8}{9.6}{\rmdefault}{\mddefault}{\updefault}{\color[rgb]{0,0,0}$z$}%
}}}}
\put(17536,-5866){\makebox(0,0)[lb]{\smash{{\SetFigFont{8}{9.6}{\rmdefault}{\mddefault}{\updefault}{\color[rgb]{0,0,0}$(z^{*})^{i+1}$}%
}}}}
\put(16816,-3706){\makebox(0,0)[lb]{\smash{{\SetFigFont{8}{9.6}{\rmdefault}{\mddefault}{\updefault}{\color[rgb]{0,0,0}$y^{*}$}%
}}}}
\put(16636,-3391){\makebox(0,0)[lb]{\smash{{\SetFigFont{8}{9.6}{\rmdefault}{\mddefault}{\updefault}{\color[rgb]{0,0,0}$y$}%
}}}}
\put(14836,-6721){\makebox(0,0)[lb]{\smash{{\SetFigFont{8}{9.6}{\rmdefault}{\mddefault}{\updefault}{\color[rgb]{0,0,0}$z^{*}$}%
}}}}
\put(16231,-6901){\makebox(0,0)[lb]{\smash{{\SetFigFont{8}{9.6}{\rmdefault}{\mddefault}{\updefault}{\color[rgb]{0,0,0}$z$}%
}}}}
\put(12991,-2581){\makebox(0,0)[lb]{\smash{{\SetFigFont{8}{9.6}{\rmdefault}{\mddefault}{\updefault}{\color[rgb]{0,0,0}$x$}%
}}}}
\put(16051,-3211){\makebox(0,0)[lb]{\smash{{\SetFigFont{8}{9.6}{\rmdefault}{\mddefault}{\updefault}{\color[rgb]{0,0,0}$y^{*}$}%
}}}}
\put(15691,-2941){\makebox(0,0)[lb]{\smash{{\SetFigFont{8}{9.6}{\rmdefault}{\mddefault}{\updefault}{\color[rgb]{0,0,0}$y$}%
}}}}
\put(16996,-3031){\makebox(0,0)[lb]{\smash{{\SetFigFont{8}{9.6}{\rmdefault}{\mddefault}{\updefault}{\color[rgb]{0,0,0}$z=xy^{*}$}%
}}}}
\put(12601,-7486){\makebox(0,0)[lb]{\smash{{\SetFigFont{8}{9.6}{\rmdefault}{\mddefault}{\updefault}{\color[rgb]{0,0,0}$A_{5}(u)=xy^{2}(y^{*})^{2}(z^{+})^{2}\cdot\mathbf{st}(u)\cdot A(u,z^{*})$}%
}}}}
\put(12151,-6001){\makebox(0,0)[lb]{\smash{{\SetFigFont{8}{9.6}{\rmdefault}{\mddefault}{\updefault}{\color[rgb]{0,0,0}$(z^{*})^{k+1}$}%
}}}}
\put(11701,-6811){\makebox(0,0)[lb]{\smash{{\SetFigFont{8}{9.6}{\rmdefault}{\mddefault}{\updefault}{\color[rgb]{0,0,0}$z$}%
}}}}
\put(11746,-3031){\makebox(0,0)[lb]{\smash{{\SetFigFont{8}{9.6}{\rmdefault}{\mddefault}{\updefault}{\color[rgb]{0,0,0}$z=xy^{*}$}%
}}}}
\put(10441,-2941){\makebox(0,0)[lb]{\smash{{\SetFigFont{8}{9.6}{\rmdefault}{\mddefault}{\updefault}{\color[rgb]{0,0,0}$y$}%
}}}}
\put(10801,-3211){\makebox(0,0)[lb]{\smash{{\SetFigFont{8}{9.6}{\rmdefault}{\mddefault}{\updefault}{\color[rgb]{0,0,0}$y^{*}$}%
}}}}
\put(7741,-2581){\makebox(0,0)[lb]{\smash{{\SetFigFont{8}{9.6}{\rmdefault}{\mddefault}{\updefault}{\color[rgb]{0,0,0}$x$}%
}}}}
\put(11971,-6496){\makebox(0,0)[lb]{\smash{{\SetFigFont{8}{9.6}{\rmdefault}{\mddefault}{\updefault}{\color[rgb]{0,0,0}$\mathbf{pi}$}%
}}}}
\put(11386,-3391){\makebox(0,0)[lb]{\smash{{\SetFigFont{8}{9.6}{\rmdefault}{\mddefault}{\updefault}{\color[rgb]{0,0,0}$y$}%
}}}}
\put(11566,-3706){\makebox(0,0)[lb]{\smash{{\SetFigFont{8}{9.6}{\rmdefault}{\mddefault}{\updefault}{\color[rgb]{0,0,0}$y^{*}$}%
}}}}
\put(8026,-7486){\makebox(0,0)[lb]{\smash{{\SetFigFont{8}{9.6}{\rmdefault}{\mddefault}{\updefault}{\color[rgb]{0,0,0}$A_{4}(u)=xy^{2}(y^{*})^{2}(z^{+})^{2}\cdot\mathbf{st}(u)\cdot\mathbf{pi}\cdot A(z^{*})$}%
}}}}
\put(10051,-7186){\makebox(0,0)[lb]{\smash{{\SetFigFont{8}{9.6}{\rmdefault}{\mddefault}{\updefault}{\color[rgb]{0,0,0}$z$}%
}}}}
\put(9301,-6586){\makebox(0,0)[lb]{\smash{{\SetFigFont{8}{9.6}{\rmdefault}{\mddefault}{\updefault}{\color[rgb]{0,0,0}$z^{*}$}%
}}}}
\put(11596,-7086){\makebox(0,0)[lb]{\smash{{\SetFigFont{8}{9.6}{\rmdefault}{\mddefault}{\updefault}{\color[rgb]{0,0,0}$\mathbf{st}(u)$}%
}}}}
\put(17076,-6431){\makebox(0,0)[lb]{\smash{{\SetFigFont{8}{9.6}{\rmdefault}{\mddefault}{\updefault}{\color[rgb]{0,0,0}$u^{k-i}$}%
}}}}
\put(16951,-6686){\makebox(0,0)[lb]{\smash{{\SetFigFont{8}{9.6}{\rmdefault}{\mddefault}{\updefault}{\color[rgb]{0,0,0}$\mathbf{st}(u)$}%
}}}}
\end{picture}%

%% file: classeB1.pstex_t
\begin{picture}(0,0)%
\includegraphics{classeB1.pstex}%
\end{picture}%
\setlength{\unitlength}{2763sp}%
\begingroup\makeatletter\ifx\SetFigFont\undefined%
\gdef\SetFigFont#1#2#3#4#5{%
  \reset@font\fontsize{#1}{#2pt}%
  \fontfamily{#3}\fontseries{#4}\fontshape{#5}%
  \selectfont}%
\fi\endgroup%
\begin{picture}(4936,5294)(8866,-7589)
\put(12286,-3481){\makebox(0,0)[lb]{\smash{{\SetFigFont{8}{9.6}{\rmdefault}{\mddefault}{\updefault}{\color[rgb]{0,0,0}$z=xy^{*}$}%
}}}}
\put(11566,-2941){\makebox(0,0)[lb]{\smash{{\SetFigFont{8}{9.6}{\rmdefault}{\mddefault}{\updefault}{\color[rgb]{0,0,0}$y$}%
}}}}
\put(11656,-3256){\makebox(0,0)[lb]{\smash{{\SetFigFont{8}{9.6}{\rmdefault}{\mddefault}{\updefault}{\color[rgb]{0,0,0}$y^{*}$}%
}}}}
\put(8866,-2581){\makebox(0,0)[lb]{\smash{{\SetFigFont{8}{9.6}{\rmdefault}{\mddefault}{\updefault}{\color[rgb]{0,0,0}$x$}%
}}}}
\put(10201,-6436){\makebox(0,0)[lb]{\smash{{\SetFigFont{8}{9.6}{\rmdefault}{\mddefault}{\updefault}{\color[rgb]{0,0,0}$z^{*}$}%
}}}}
\put(12796,-7052){\makebox(0,0)[lb]{\smash{{\SetFigFont{8}{9.6}{\rmdefault}{\mddefault}{\updefault}{\color[rgb]{0,0,0}$z$}%
}}}}
\put(13156,-4186){\makebox(0,0)[lb]{\smash{{\SetFigFont{8}{9.6}{\rmdefault}{\mddefault}{\updefault}{\color[rgb]{0,0,0}$y^{*}$}%
}}}}
\put(9091,-7531){\makebox(0,0)[lb]{\smash{{\SetFigFont{8}{9.6}{\rmdefault}{\mddefault}{\updefault}{\color[rgb]{0,0,0}$B_1(u)=xy^+y^*z^+\cdot yuz^*\cdot\mathbf{st}(u)\cdot z\cdot\mathbf{pi}\cdot A(u)$}%
}}}}
\put(10471,-7141){\makebox(0,0)[lb]{\smash{{\SetFigFont{8}{9.6}{\rmdefault}{\mddefault}{\updefault}{\color[rgb]{0,0,0}$z\cdot\mathbf{pi}$}%
}}}}
\put(13261,-5896){\makebox(0,0)[lb]{\smash{{\SetFigFont{8}{9.6}{\rmdefault}{\mddefault}{\updefault}{\color[rgb]{0,0,0}$u^k$}%
}}}}
\put(13166,-6421){\makebox(0,0)[lb]{\smash{{\SetFigFont{8}{9.6}{\rmdefault}{\mddefault}{\updefault}{\color[rgb]{0,0,0}$yuz^*\cdot\mathbf {st}(u)$}%
}}}}
\end{picture}%

%% file: classeB4.pstex_t
\begin{picture}(0,0)%
\includegraphics{classeB4.pstex}%
\end{picture}%
\setlength{\unitlength}{2763sp}%
\begingroup\makeatletter\ifx\SetFigFont\undefined%
\gdef\SetFigFont#1#2#3#4#5{%
  \reset@font\fontsize{#1}{#2pt}%
  \fontfamily{#3}\fontseries{#4}\fontshape{#5}%
  \selectfont}%
\fi\endgroup%
\begin{picture}(5518,5034)(9886,-7335)
\put(9886,-2467){\makebox(0,0)[lb]{\smash{{\SetFigFont{8}{9.6}{\rmdefault}{\mddefault}{\updefault}{\color[rgb]{0,0,0}$x$}%
}}}}
\put(14026,-2986){\makebox(0,0)[lb]{\smash{{\SetFigFont{8}{9.6}{\rmdefault}{\mddefault}{\updefault}{\color[rgb]{0,0,0}$z=xy^{*}$}%
}}}}
\put(13801,-3436){\makebox(0,0)[lb]{\smash{{\SetFigFont{8}{9.6}{\rmdefault}{\mddefault}{\updefault}{\color[rgb]{0,0,0}$y$}%
}}}}
\put(13801,-3736){\makebox(0,0)[lb]{\smash{{\SetFigFont{8}{9.6}{\rmdefault}{\mddefault}{\updefault}{\color[rgb]{0,0,0}$y^{*}$}%
}}}}
\put(12451,-2986){\makebox(0,0)[lb]{\smash{{\SetFigFont{8}{9.6}{\rmdefault}{\mddefault}{\updefault}{\color[rgb]{0,0,0}$y$}%
}}}}
\put(12526,-3211){\makebox(0,0)[lb]{\smash{{\SetFigFont{8}{9.6}{\rmdefault}{\mddefault}{\updefault}{\color[rgb]{0,0,0}$y^{*}$}%
}}}}
\put(14081,-5782){\makebox(0,0)[lb]{\smash{{\SetFigFont{8}{9.6}{\rmdefault}{\mddefault}{\updefault}{\color[rgb]{0,0,0}$(z^{*})^{i+1}$}%
}}}}
\put(13871,-6142){\makebox(0,0)[lb]{\smash{{\SetFigFont{8}{9.6}{\rmdefault}{\mddefault}{\updefault}{\color[rgb]{0,0,0}$z$}%
}}}}
\put(14161,-6452){\makebox(0,0)[lb]{\smash{{\SetFigFont{8}{9.6}{\rmdefault}{\mddefault}{\updefault}{\color[rgb]{0,0,0}$u^{k-i-j}$}%
}}}}
\put(13866,-6642){\makebox(0,0)[lb]{\smash{{\SetFigFont{8}{9.6}{\rmdefault}{\mddefault}{\updefault}{\color[rgb]{0,0,0}$z$}%
}}}}
\put(13106,-6872){\makebox(0,0)[lb]{\smash{{\SetFigFont{8}{9.6}{\rmdefault}{\mddefault}{\updefault}{\color[rgb]{0,0,0}$(z^{*})^{j+1}$, $j\neq0$}%
}}}}
\put(10131,-7277){\makebox(0,0)[lb]{\smash{{\SetFigFont{8}{9.6}{\rmdefault}{\mddefault}{\updefault}{\color[rgb]{0,0,0}$B_{4}(u)=x(y^{+})^2(1+z\cdot\mathbf{pi})(z^+)^2(A(z^*,z^*,u)-A(z^*,u))$}%
}}}}
\put(10446,-6542){\makebox(0,0)[lb]{\smash{{\SetFigFont{8}{9.6}{\rmdefault}{\mddefault}{\updefault}{\color[rgb]{0,0,0}$1+z\cdot\mathbf{pi}$}%
}}}}
\end{picture}%

%% file: classeB5-7.pstex_t
\begin{picture}(0,0)%
\includegraphics{classeB5-7.pstex}%
\end{picture}%
\setlength{\unitlength}{2763sp}%
\begingroup\makeatletter\ifx\SetFigFont\undefined%
\gdef\SetFigFont#1#2#3#4#5{%
  \reset@font\fontsize{#1}{#2pt}%
  \fontfamily{#3}\fontseries{#4}\fontshape{#5}%
  \selectfont}%
\fi\endgroup%
\begin{picture}(11833,4965)(3591,-7260)
\put(9106,-2957){\makebox(0,0)[lb]{\smash{{\SetFigFont{8}{9.6}{\rmdefault}{\mddefault}{\updefault}{\color[rgb]{0,0,0}$y$}%
}}}}
\put(6616,-2477){\makebox(0,0)[lb]{\smash{{\SetFigFont{8}{9.6}{\rmdefault}{\mddefault}{\updefault}{\color[rgb]{0,0,0}$x$}%
}}}}
\put(9841,-6002){\makebox(0,0)[lb]{\smash{{\SetFigFont{8}{9.6}{\rmdefault}{\mddefault}{\updefault}{\color[rgb]{0,0,0}$(z^{*})^{i+1}$}%
}}}}
\put(6766,-7202){\makebox(0,0)[lb]{\smash{{\SetFigFont{8}{9.6}{\rmdefault}{\mddefault}{\updefault}{\color[rgb]{0,0,0}$B_{6}(u)=2\cdot xy^*y^{+}z^*B(z^*,u)$}%
}}}}
\put(9766,-3932){\makebox(0,0)[lb]{\smash{{\SetFigFont{8}{9.6}{\rmdefault}{\mddefault}{\updefault}{\color[rgb]{0,0,0}$z=xy^{*}$}%
}}}}
\put(9181,-3242){\makebox(0,0)[lb]{\smash{{\SetFigFont{8}{9.6}{\rmdefault}{\mddefault}{\updefault}{\color[rgb]{0,0,0}$y^{*}$}%
}}}}
\put(9256,-3617){\makebox(0,0)[lb]{\smash{{\SetFigFont{8}{9.6}{\rmdefault}{\mddefault}{\updefault}{\color[rgb]{0,0,0}$y^{*}$}%
}}}}
\put(9316,-6332){\makebox(0,0)[lb]{\smash{{\SetFigFont{8}{9.6}{\rmdefault}{\mddefault}{\updefault}{\color[rgb]{0,0,0}$u^{k-i}$}%
}}}}
\put(5701,-4036){\makebox(0,0)[lb]{\smash{{\SetFigFont{8}{9.6}{\rmdefault}{\mddefault}{\updefault}{\color[rgb]{0,0,0}$y^{*}$}%
}}}}
\put(5656,-3436){\makebox(0,0)[lb]{\smash{{\SetFigFont{8}{9.6}{\rmdefault}{\mddefault}{\updefault}{\color[rgb]{0,0,0}$y^{*}$}%
}}}}
\put(3901,-7186){\makebox(0,0)[lb]{\smash{{\SetFigFont{8}{9.6}{\rmdefault}{\mddefault}{\updefault}{\color[rgb]{0,0,0}$B_{5}(u)=x(y^{*})^2B(u)$}%
}}}}
\put(3601,-2461){\makebox(0,0)[lb]{\smash{{\SetFigFont{8}{9.6}{\rmdefault}{\mddefault}{\updefault}{\color[rgb]{0,0,0}$x$}%
}}}}
\put(10681,-2477){\makebox(0,0)[lb]{\smash{{\SetFigFont{8}{9.6}{\rmdefault}{\mddefault}{\updefault}{\color[rgb]{0,0,0}$x$}%
}}}}
\put(10831,-7202){\makebox(0,0)[lb]{\smash{{\SetFigFont{8}{9.6}{\rmdefault}{\mddefault}{\updefault}{\color[rgb]{0,0,0}$B_{7}(u)=x(y^{+})^2(z^+)^2B(z^*,z^*,u)$}%
}}}}
\put(13246,-3242){\makebox(0,0)[lb]{\smash{{\SetFigFont{8}{9.6}{\rmdefault}{\mddefault}{\updefault}{\color[rgb]{0,0,0}$y^{*}$}%
}}}}
\put(13876,-3722){\makebox(0,0)[lb]{\smash{{\SetFigFont{8}{9.6}{\rmdefault}{\mddefault}{\updefault}{\color[rgb]{0,0,0}$y^{*}$}%
}}}}
\put(13816,-3452){\makebox(0,0)[lb]{\smash{{\SetFigFont{8}{9.6}{\rmdefault}{\mddefault}{\updefault}{\color[rgb]{0,0,0}$y$}%
}}}}
\put(13186,-2987){\makebox(0,0)[lb]{\smash{{\SetFigFont{8}{9.6}{\rmdefault}{\mddefault}{\updefault}{\color[rgb]{0,0,0}$y$}%
}}}}
\put(13711,-2942){\makebox(0,0)[lb]{\smash{{\SetFigFont{8}{9.6}{\rmdefault}{\mddefault}{\updefault}{\color[rgb]{0,0,0}$z=xy^{*}$}%
}}}}
\put(14101,-5747){\makebox(0,0)[lb]{\smash{{\SetFigFont{8}{9.6}{\rmdefault}{\mddefault}{\updefault}{\color[rgb]{0,0,0}$(z^{*})^{i+1}$}%
}}}}
\put(13006,-6857){\makebox(0,0)[lb]{\smash{{\SetFigFont{8}{9.6}{\rmdefault}{\mddefault}{\updefault}{\color[rgb]{0,0,0}$(z^{*})^{j+1}$}%
}}}}
\put(13891,-6107){\makebox(0,0)[lb]{\smash{{\SetFigFont{8}{9.6}{\rmdefault}{\mddefault}{\updefault}{\color[rgb]{0,0,0}$z$}%
}}}}
\put(14101,-6377){\makebox(0,0)[lb]{\smash{{\SetFigFont{8}{9.6}{\rmdefault}{\mddefault}{\updefault}{\color[rgb]{0,0,0}$u^{k-i-j}$}%
}}}}
\put(13906,-6647){\makebox(0,0)[lb]{\smash{{\SetFigFont{8}{9.6}{\rmdefault}{\mddefault}{\updefault}{\color[rgb]{0,0,0}$z$}%
}}}}
\put(5671,-5912){\makebox(0,0)[lb]{\smash{{\SetFigFont{8}{9.6}{\rmdefault}{\mddefault}{\updefault}{\color[rgb]{0,0,0}$u^{k}$}%
}}}}
\end{picture}%

%% file: classeP1-3.pstex_t
\begin{picture}(0,0)%
\includegraphics{classeP1-3.pstex}%
\end{picture}%
\setlength{\unitlength}{2763sp}%
\begingroup\makeatletter\ifx\SetFigFont\undefined%
\gdef\SetFigFont#1#2#3#4#5{%
  \reset@font\fontsize{#1}{#2pt}%
  \fontfamily{#3}\fontseries{#4}\fontshape{#5}%
  \selectfont}%
\fi\endgroup%
\begin{picture}(11658,4649)(1491,-6944)
\put(3556,-3436){\makebox(0,0)[lb]{\smash{{\SetFigFont{8}{9.6}{\rmdefault}{\mddefault}{\updefault}{\color[rgb]{0,0,0}$y^{*}$}%
}}}}
\put(1576,-6286){\makebox(0,0)[lb]{\smash{{\SetFigFont{8}{9.6}{\rmdefault}{\mddefault}{\updefault}{\color[rgb]{0,0,0}$P_{1}=xy^{*}A(1)$}%
}}}}
\put(1501,-2911){\makebox(0,0)[lb]{\smash{{\SetFigFont{8}{9.6}{\rmdefault}{\mddefault}{\updefault}{\color[rgb]{0,0,0}$x$}%
}}}}
\put(10441,-2941){\makebox(0,0)[lb]{\smash{{\SetFigFont{8}{9.6}{\rmdefault}{\mddefault}{\updefault}{\color[rgb]{0,0,0}$y$}%
}}}}
\put(10801,-3211){\makebox(0,0)[lb]{\smash{{\SetFigFont{8}{9.6}{\rmdefault}{\mddefault}{\updefault}{\color[rgb]{0,0,0}$y^{*}$}%
}}}}
\put(7741,-2581){\makebox(0,0)[lb]{\smash{{\SetFigFont{8}{9.6}{\rmdefault}{\mddefault}{\updefault}{\color[rgb]{0,0,0}$x$}%
}}}}
\put(11776,-3286){\makebox(0,0)[lb]{\smash{{\SetFigFont{8}{9.6}{\rmdefault}{\mddefault}{\updefault}{\color[rgb]{0,0,0}$z=xy^{*}$}%
}}}}
\put(6841,-2941){\makebox(0,0)[lb]{\smash{{\SetFigFont{8}{9.6}{\rmdefault}{\mddefault}{\updefault}{\color[rgb]{0,0,0}$y$}%
}}}}
\put(6931,-3256){\makebox(0,0)[lb]{\smash{{\SetFigFont{8}{9.6}{\rmdefault}{\mddefault}{\updefault}{\color[rgb]{0,0,0}$y^{*}$}%
}}}}
\put(4141,-2581){\makebox(0,0)[lb]{\smash{{\SetFigFont{8}{9.6}{\rmdefault}{\mddefault}{\updefault}{\color[rgb]{0,0,0}$x$}%
}}}}
\put(7951,-6886){\makebox(0,0)[lb]{\smash{{\SetFigFont{8}{9.6}{\rmdefault}{\mddefault}{\updefault}{\color[rgb]{0,0,0}$P_3=xy^{+}(1+z\cdot\mathbf{pi})(z^{*}A(z^{*})-A(1))$}%
}}}}
\put(7126,-3586){\makebox(0,0)[lb]{\smash{{\SetFigFont{8}{9.6}{\rmdefault}{\mddefault}{\updefault}{\color[rgb]{0,0,0}$z=xy^{*}$}%
}}}}
\put(11626,-6436){\makebox(0,0)[lb]{\smash{{\SetFigFont{8}{9.6}{\rmdefault}{\mddefault}{\updefault}{\color[rgb]{0,0,0}$(z^{*})^{k+1}-1$}%
}}}}
\put(8581,-6436){\makebox(0,0)[lb]{\smash{{\SetFigFont{8}{9.6}{\rmdefault}{\mddefault}{\updefault}{\color[rgb]{0,0,0}$1+z\cdot\mathbf{pi}$}%
}}}}
\put(5601,-6436){\makebox(0,0)[lb]{\smash{{\SetFigFont{8}{9.6}{\rmdefault}{\mddefault}{\updefault}{\color[rgb]{0,0,0}$z\cdot\mathbf {pi}$}%
}}}}
\put(4276,-6886){\makebox(0,0)[lb]{\smash{{\SetFigFont{8}{9.6}{\rmdefault}{\mddefault}{\updefault}{\color[rgb]{0,0,0}$P_{2}(u)=xy^{+}z\cdot\mathbf{pi}\cdot A(1)$}%
}}}}
\end{picture}%

%% file: classeP4-5.pstex_t
\begin{picture}(0,0)%
\includegraphics{classeP4-5.pstex}%
\end{picture}%
\setlength{\unitlength}{2763sp}%
\begingroup\makeatletter\ifx\SetFigFont\undefined%
\gdef\SetFigFont#1#2#3#4#5{%
  \reset@font\fontsize{#1}{#2pt}%
  \fontfamily{#3}\fontseries{#4}\fontshape{#5}%
  \selectfont}%
\fi\endgroup%
\begin{picture}(9108,4949)(3591,-7244)
\put(10441,-2941){\makebox(0,0)[lb]{\smash{{\SetFigFont{8}{9.6}{\rmdefault}{\mddefault}{\updefault}{\color[rgb]{0,0,0}$y$}%
}}}}
\put(5656,-3436){\makebox(0,0)[lb]{\smash{{\SetFigFont{8}{9.6}{\rmdefault}{\mddefault}{\updefault}{\color[rgb]{0,0,0}$y^{*}$}%
}}}}
\put(3901,-7186){\makebox(0,0)[lb]{\smash{{\SetFigFont{8}{9.6}{\rmdefault}{\mddefault}{\updefault}{\color[rgb]{0,0,0}$P_{4}=xy^{*}B(1)$}%
}}}}
\put(7951,-2461){\makebox(0,0)[lb]{\smash{{\SetFigFont{8}{9.6}{\rmdefault}{\mddefault}{\updefault}{\color[rgb]{0,0,0}$x$}%
}}}}
\put(3601,-2461){\makebox(0,0)[lb]{\smash{{\SetFigFont{8}{9.6}{\rmdefault}{\mddefault}{\updefault}{\color[rgb]{0,0,0}$x$}%
}}}}
\put(10576,-3286){\makebox(0,0)[lb]{\smash{{\SetFigFont{8}{9.6}{\rmdefault}{\mddefault}{\updefault}{\color[rgb]{0,0,0}$y^{*}$}%
}}}}
\put(11101,-3736){\makebox(0,0)[lb]{\smash{{\SetFigFont{8}{9.6}{\rmdefault}{\mddefault}{\updefault}{\color[rgb]{0,0,0}$z=xy^{*}$}%
}}}}
\put(11176,-5986){\makebox(0,0)[lb]{\smash{{\SetFigFont{8}{9.6}{\rmdefault}{\mddefault}{\updefault}{\color[rgb]{0,0,0}$(z^{*})^{k+1}-1$}%
}}}}
\put(8101,-7186){\makebox(0,0)[lb]{\smash{{\SetFigFont{8}{9.6}{\rmdefault}{\mddefault}{\updefault}{\color[rgb]{0,0,0}$P_{5}=xy^{+}(z^*B(z^*)-B(1))$}%
}}}}
\end{picture}%